\documentclass[twoside,12pt]{article}

\usepackage{amssymb,amsmath}

\allowdisplaybreaks[1]

\newenvironment{proof}{\begin{trivlist}\item[]{\it
Proof.}}{\hfill$\square$\end{trivlist}}

\newtheorem{theorem}{Theorem}[section]

\newtheorem{proposition}[theorem]{Proposition}
\newtheorem{remark}[theorem]{Remark}

\def\mc{{\mathbb{C}}}
\def\mz{{\mathbb{Z}}}
\def\qtr{{\mathrm{Tr}}_q}
\def\Tr{{\mathrm{Tr}}}
\def\qdet{\det_q}
\def\AA{{\mathcal{A}}}
\def\genmat{{\mathcal{G}}(N,m)}
\def\qgenmat{{\mathcal{G}}_q(N,m)}
\def\twotwogenmat{{\mathcal{G}}(2,2)}
\def\twotwoqgenmat{{\mathcal{G}}_q(2,2)}
\def\matinv{{\mathcal{R}}(N,m)}
\def\twotwomatinv{{\mathcal{R}}(2,2)}
\def\qmatinv{{\mathcal{R}}_q(N,m)}
\def\twomqmatinv{{\mathcal{R}}_q(2,m)}
\def\twotwoqmatinv{{\mathcal{R}}_q(2,2)}
\def\tracering{{\mathcal{T}}(N,m)}
\def\twotwotracering{{\mathcal{T}}(2,2)}
\def\qtracering{{\mathcal{T}}_q(N,m)}
\def\twotwoqtracering{{\mathcal{T}}_q(2,2)}
\def\twomqtracering{{\mathcal{T}}_q(2,m)}
\def\qsl{{\mathcal{F}}_q(SL_N)}
\def\twoqsl{{\mathcal{F}}_q(SL_2)}
\def\funsl{{\mathcal{F}}(SL_N)}
\def\sevmat{{\mathcal{A}}(N,m)}
\def\qsevmat{{\mathcal{A}}_q(N,m)}
\def\ntwoqsevmat{{\mathcal{A}}_q(N,2)}
\def\nthreeqsevmat{{\mathcal{A}}_q(N,3)}
\def\twomqsevmat{{\mathcal{A}}_q(2,m)}
\def\twothreeqsevmat{{\mathcal{A}}_q(2,3)}
\def\twotwoqsevmat{{\mathcal{A}}_q(2,2)}
\def\mcoact{\beta(m)}
\def\coact{\beta}
\def\funmn{{\mathcal{F}}(M_N)}
\def\qfunmn{{\mathcal{F}}_q(M_N)}
\def\rea{{\mathcal{A}}_q(N)}
\def\tworea{{\mathcal{A}}_q(2)}
\def\btensor{{\underline{\otimes}}}
\def\coaction{\beta}
\def\rhat{\widehat{R}}
\def\rtilde{\widetilde{R}}
\def\fundiag{{\mathcal{F}}(K)}
\def\id{{\mathrm{id}}}
\def\xijspan{{\mathcal{X}}}
\def\linspan{{\mathrm{Span}}_{\mc}}
\def\triv{{\mathrm{triv}}}
\def\Endom{{\mathrm{End}}}
\def\Homom{{\mathrm{Hom}}}
\begin{document}
\date{}
\title{Quantized trace rings}
\author{M. Domokos 
\thanks{Supported through a European Community Marie Curie Fellowship. 
Partially supported by OTKA No. T034530, T046378, and 
a Bolyai Research Fellowship.}\;
and T. H. Lenagan 
\thanks{Both authors partially supported by Leverhulme Research 
Interchange Grant F/00158/X.}}
\maketitle

\begin{abstract} The general linear group acts on $m$-tuples of 
$N\times N$ matrices by simultaneous conjugation. 
Quantum deformations of the corresponding rings of invariants 
and the so-called trace rings are investigated. 
\end{abstract}

\medskip
\noindent 2000 Mathematics Subject Classification. 16W35, 16R30, 20G42


\section{Introduction}\label{sec:intro} 

Write $M_N^m$ for the space 
$M_N(\mc)\oplus\cdots\oplus M_N(\mc)$ 
of $m$-tuples of $N\times N$ complex matrices. 
Denote by $U(k)=\left(u(k)^i_j\right)_{i,j=1}^N$ 
the function mapping an $m$-tuple to its $k$th component $(k=1,\ldots,m)$. 
The $N^2m$-variable commutative polynomial algebra 
$\AA=\sevmat$ 
generated by the $u(k)^i_j$ is the coordinate ring of $M_N^m$. 
Let us recall various subalgebras of $\AA$ and $M_N(\AA)$. 
The $\mc$-subalgebra of $M_N(\AA)$ generated by $U(1),\ldots,U(m)$ and 
the identity matrix, $I$, is known as the algebra 
$\genmat$ of $m$ generic $N\times N$ 
matrices. The complex special linear group $SL_N$ 
acts on $M_N^m$ by simultaneous conjugation. This induces an action 
on $\AA$ by linear substitution of the variables: 
$g\in SL_N$ maps $u(k)^i_j$ to the $(i,j)$-entry of $g^{-1}U(k)g$. 
The ring $\matinv$ of {\it matrix invariants} consists of 
\[\matinv=
\{f\in\sevmat \mid \forall g\in SL_N:\ g\cdot f=f\}.\] 
Let $\tracering$ denote the algebra of $SL_N$-equivariant polynomial maps from 
$M_N^m$ to $M_N(\mc)$. Identify polynomial maps 
$M_N^m\to M_N(\mc)$ with elements of $M_N(\AA)$ in the usual way.  
Then $\tracering$ is the $\matinv$-subalgebra of $M_N(\AA)$ 
generated by $I,U(1),\ldots,U(m)$. It is called the 
{\it trace ring} (or sometimes the {\it algebra of matrix concomitants}). 
The reason for the name is that $\matinv$ is generated by elements of the form 
$\Tr(U(k_1)\cdots U(k_d))$, where 
$k_1,\ldots,k_d\in\{1,\ldots,m\}$, $d\leq N^2$, and $\Tr$ 
denotes the usual trace function (see \cite{sib}, \cite{raz}, \cite{pr}). 

The algebra $\genmat$ of generic matrices appeared in the theory of algebras 
satisfying polynomial identities, as a relatively free object. 
On the other hand, $\matinv$ and $\tracering$ are natural objects 
of invariant theory. By definition, $\genmat\subset \tracering$; 
although this is a proper inclusion, $\genmat$ turns out to be closely 
related to $\tracering$. As a remarkable consequence, certain questions in 
non-commutative algebra can be approached by methods of classical invariant 
theory, see \cite{pr} and \cite{fo} for details. 

Our aim in the present paper is to find analogues of the algebras 
$\matinv$, $\tracering$, $\genmat$ in the context of quantum groups. 
Write $\qsl$ for the quantized function algebra of $SL_N$, 
defined for example in \cite{rtf}. Here $q\in\mc^{\times}$ 
is a non-zero parameter. 
$\qsl$ is a Hopf algebra, which is a non-commutative deformation 
of the coordinate ring  $\funsl$ of $SL_N$. 
(For unexplained terminology and facts relating to Hopf algebras, 
coactions, and quantum groups, 
we refer to the books \cite{kl-sch}, \cite{majid2}.)   
Works of Majid \cite{majid2} and 
Kulish and Sasaki \cite{ku-sa} provide us with an appropriate 
non-commutative deformation $\qsevmat$ of $\sevmat$, together with 
a right coaction 
$\mcoact:\qsevmat\to\qsevmat\otimes\qsl$, 
which is a deformation of the comorphism 
$\sevmat\to\sevmat\otimes\funsl$ 
of the simultaneous conjugation action of $SL_N$ on $M_N^m$. 
To simplify notation, when $m=1$ we write $\rea$ instead of 
${\mathcal{A}}_q(N,1)$. 
We take for the algebra $\rea$ the {\it braided matrices} of 
\cite{majid1} (called the {\it reflection equation algebra} in 
\cite{ku-skl}), associated with the R-matrix belonging to 
the quantum $SL_N$. 
It is a right $\qsl$-comodule algebra, and for some problems relating the 
conjugation action of $SL_N$ on $M_N(\mc)$, 
it appears to be a better quantum deformation of $\funmn$, 
than the coordinate ring $\qfunmn$ of quantum matrices, 
see for example \cite{gu-sa}, \cite{do-mu}. 
When $m=2$, we take for $\ntwoqsevmat$ the {\it braided tensor product} 
(cf. \cite{majid2}) 
$\rea\btensor\rea$ 
of two copies of $\rea$. 
As a vector space, it is just the ordinary tensor product 
$\rea\otimes_{\mc}\rea$, 
endowed with the tensor product $\qsl$-comodule structure. 
Using that $\qsl$ is coquasitriangular, 
one can define a multiplication on this tensor product space 
that makes the coaction multiplicative, see 
Corollary 9.2.14 in \cite{majid2}, and the sign $\btensor$ indicates 
the vector space tensor product endowed with this multiplication. 
With an iterated use of braided tensor products we arrive at $\qsevmat$ 
for arbitrary $m$. The same algebra was defined in \cite{ku-sa} 
in terms of generators and relations. 

Now we place the algebra $\qmatinv$ of $\mcoact$-coinvariants into the 
centre of 
investigation, as well as a subalgebra $\qtracering$ of $M_N(\qsevmat)$.   
Under certain assumptions on $q$ we show that 
these algebras can be considered as quantum versions of the ring of matrix 
invariants and the trace ring, by investigating their Hilbert series, 
their generators, and ring theoretical properties. 
In the special case $N=2$ we present some explicit calculations, 
leading to the somewhat surprising conclusion that 
$\twotwoqtracering\cong\twotwotracering$, and 
$\twotwoqgenmat\cong\twotwogenmat$. 


\section{Some quadratic algebras}\label{sec:quadratic} 

To define $\qsl$ one starts with the R-matrix 
$R=\left(R^{ik}_{jl}\right)_{i,j,k,l=1}^N$. 
This is an $N^2\times N^2$ matrix, whose rows and columns are indexed by 
pairs, and the entry in the crossing of the $(i,k)$ row and the 
$(j,l)$ column is 
\[R^{ik}_{jl}=\begin{cases} 
q, &\mbox{ if }i=j=k=l;\\
1, &\mbox{ if }i=j,k=l,i\neq k;\\
q-q^{-1},&\mbox{ if }i>j,i=l,j=k; \\
0, &\mbox{ otherwise.} 
\end{cases}\]
It is usual to think of $R$ as the matrix of a linear endomorphism of 
$V\otimes V$ with respect to the basis 
$e_1\otimes e_1,e_1\otimes e_2,\ldots,e_1\otimes e_N,e_2\otimes e_1,\ldots,
e_N\otimes e_N$. 
So identifying $R$ and the endomorphism we have 
$R(e_j\otimes e_l)=R^{ik}_{jl}e_i\otimes e_k$, 
where the convention of summing over repeated indices is used 
(and will be used throughout the paper). 
 
The {\it quantized coordinate ring} $\qfunmn$ is the $\mc$-algebra 
generated by the entries of $T=\left(t^i_j\right)_{i,j=1}^N$, 
subject to the relations 
$RT_1T_2=T_2T_1R$, 
with the matrix notation of \cite{rtf}; that is, $T_1$, $T_2$ are 
the matrix Kronecker products 
$T_1=T\otimes I$, $T_2=I\otimes T$, 
so 
$(T_1)^{ik}_{jl}=t^i_j\delta^k_l$, 
$(T_2)^{ik}_{jl}=\delta^i_jt^k_l$, 
where $\delta^i_i=1$, and $\delta^i_j=0$ for $i\neq j$. 
The algebra $\qfunmn$ is a bialgebra with comultiplication 
$\Delta(t^i_j)=t^i_k\otimes t^k_j$ 
and counit $\varepsilon(t^i_j)=\delta^i_j$. 
The quantum determinant 
$\qdet=\sum_{\pi\in S_N}(-q)^{{\mathrm{length}}(\pi)}
t^1_{\pi(1)}\cdots t^N_{\pi(N)}$ 
is a central group-like element, and 
$\qsl$ is defined as the quotient of $\qfunmn$ modulo the ideal 
$\langle \qdet-1\rangle$. We keep the notation $t^i_j\in\qsl$ 
for the images of the generators under the natural surjection. 
Then $\qsl$ is a Hopf algebra with antipode $S$ satisfying 
$S(T)T=I=TS(T)$. 

Denote by $\tau$ the flip endomorphism $\tau(v\otimes w)=w\otimes v$ 
of $V\otimes V$.  
Define 
$\widehat{R}=\tau\circ R$ and $R_{21}=\tau\circ R\circ\tau$. 
In coordinates, we have 
$\widehat{R}^{ik}_{jl}=R^{ki}_{jl}$ and 
$(R_{21})^{ik}_{jl}=R^{ki}_{lj}$. 
The {\it reflection equation algebra} 
$\rea$ is the $\mc$-algebra generated by the entries of 
$X=\left(x^i_j\right)_{i,j=1}^N$ 
subject to the quadratic relations 
\begin{equation} \label{eq:rea}
X_2\widehat{R}X_2\widehat{R}=\widehat{R}X_2\widehat{R}X_2. 
\end{equation}
The relations can be written also in the equivalent form 
\[X_2R_{21}X_1R=R_{21}X_1RX_2.\] 
By Proposition 10.3.5 in \cite{majid2}, 
$\rea$ is a right $\qsl$-comodule algebra in a natural way. 
Namely, the map 
$X\mapsto S(T)XT$
extends to an algebra homomorphism 
$\coaction:\rea\to\rea\otimes\qsl$, 
and $\coaction$ is a right coaction of the Hopf algebra 
$\qsl$ on $\rea$. 
Here we treat the entries of $S(T)XT$ 
(and more generally, products of the form $txu$ with 
$t,u\in\qsl$, $x\in\rea$) 
as elements of $\rea\otimes\qsl$ 
by assuming that the generators of $\qsl$ commute with the generators of 
$\rea$, and suppressing the $\otimes$ sign from the notation; so  
$\coact(x^i_j)=x^k_l \otimes S(t^i_k)t^l_j$.  
Taking into account that $S(T)=T^{-1}$, this notational device 
will be very convenient later. 

Next take two copies of $\rea$. The matrix generators of them will be 
denoted by $X=\left(x^i_j\right)$ and $Y=\left(y^i_j\right)$, respectively. 
Define $\ntwoqsevmat$ as the algebra generated by the 
$2N^2$ elements $x^i_j$, $y^i_j$, subject to the additional relations 
\begin{equation}\label{eq:qsevmatntwo} 
R^{-1}Y_1RX_2=X_2R^{-1}Y_1R.
\end{equation} 
By Theorem 10.3.1 of \cite{majid2} this algebra can be identified with the 
braided tensor product 
$\rea\btensor\rea$, 
by setting 
$x^i_j\mapsto x^i_j\btensor 1$, 
$y^i_j\mapsto 1\btensor x^i_j$. 
The relations \eqref{eq:qsevmatntwo} are equivalent to 
\begin{equation}\label{eq:qsevmatntwo-coord} 
\forall i,j,k,l:\quad 
y^i_jx^k_l=
(R^{-1})^{an}_{rc}R^{sc}_{bl}R^{id}_{am}(\rtilde)^{bk}_{jd}
x^m_ny^r_s,
\end{equation} 
where 
$\rtilde^{t_2}=(R^{t_2})^{-1}$,  
and for an $N^2\times N^2$ matrix $A$, the matrix $A^{t_2}$ 
is defined by 
$(A^{t_2})^{ik}_{jl}=A^{il}_{jk}$. 

More generally, take $m$ copies of $\rea$ with matrix generators 
\[X(1)=\left(x(1)^i_j\right),\ldots,X(m)=\left(x(m)^i_j\right),\] 
respectively, and define $\qsevmat$ to be the algebra generated 
by $\{x(k)^i_j\mid k=1,\ldots,m; 
\ i,j=1,\ldots,N\}$, 
where in addition to the reflection equation 
\eqref{eq:rea} imposed for each $X(k)=X$, we impose the relations  
\begin{equation}\label{eq:qsevmat} 
\forall r<s: \quad
R^{-1}X(s)_1RX(r)_2=X(r)_2R^{-1}X(s)_1R. 
\end{equation} 
In other words, for each $r<s$ the subalgebra of $\qsevmat$ generated by the 
entries of $X(r)$ and $X(s)$ is isomorphic to $\ntwoqsevmat$ via an isomorphism 
$x(r)^i_j\mapsto x^i_j$, $x(s)^i_j\mapsto y^i_j$. 
This algebra was defined in \cite{ku-sa}. 
It is clearly the same as the braided tensor product 
$\rea\btensor\cdots\btensor\rea$ ($m$ copies) in the sense of 
\cite{majid2}. 
Hence we may define the coaction 
$\mcoact:\qsevmat\to\qsevmat\otimes\qsl$ 
by taking the $m$th tensor power of $\coaction$. 
The multiplication in the braided tensor product is designed to make 
$\mcoact$ an algebra homomorphism. Therefore $\qsevmat$ becomes an 
$\qsl$-comodule algebra. 
More explicitly, the coaction is given on the generators by 
$\mcoact:X(r)\mapsto S(T)X(r)T$ for 
$r=1,\ldots,m$. 

\section{Polynormal sequences in the reflection equation algebra}
\label{seq:polynormal} 

It is standard that the $\mz$-graded algebra $\qfunmn$ has the 
same Hilbert series as $\funmn$. 
On the other hand, a vector space isomorphism between the degree $d$ 
homogeneous components of $\qfunmn$ and $\rea$ is established 
in (7.37) of \cite{majid2}. Consequently, the graded algebra $\rea$ 
has the same Hilbert series as the $N^2$-variable commutative polynomial 
algebra. 
We refine this statement by having a closer look at the defining relations 
of $\rea$. 
The matrix equality \eqref{eq:rea} is a shorthand for the following: 
\begin{equation} \label{eq:rea-coord} 
\forall i,j,k,l:\quad 
\rhat^{ib}_{st}\rhat^{sa}_{jl}x^k_bx^t_a
=\rhat^{ik}_{sa}\rhat^{st}_{jb}x^a_tx^b_l 
\end{equation} 
Taking into account the exact form of the R-matrix, 
one gets the following list 
of relations. 
\begin{eqnarray*}\label{eq:rea-detail} 
x^i_jx^i_l&=&qx^i_lx^i_j \qquad\mbox{ for }j<l,i\neq j,i\neq l;\\  
q^{-1}x^i_jx^i_i-qx^i_ix^i_j
&=&(q-q^{-1})\sum_{s>i}x^i_sx^s_j \qquad\mbox{ for }i>j;\\  
x^i_lx^i_i-x^i_ix^i_l&=&(1-q^{-2})\sum_{s>i}x^i_sx^s_l 
\qquad\mbox{ for }i<l;\\  
x^k_jx^i_j&=&qx^i_jx^k_j \qquad\mbox{ for }i<k,j\neq i,j\neq k;\\  
x^i_ix^k_i-x^k_ix^i_i&=&(1-q^{-2})\sum_{s>i}x^k_sx^s_i
\qquad\mbox{ for }i<k;\\  
q^{-1}x^j_jx^i_j-qx^i_jx^j_j&=&(q-q^{-1})\sum_{s>j}x^i_sx^s_j 
\qquad\mbox{ for }i<j;\\ 
x^i_jx^k_l&=&x^k_lx^i_j \qquad\mbox{ for }i<k,j<l,i\neq l,j\neq k;\\  
x^j_lx^i_j-qx^i_jx^j_l&=&(q-q^{-1})\sum_{s>j}x^i_sx^s_l 
\qquad\mbox{ for }i<j<l;\\  
x^i_jx^k_i-qx^k_ix^i_j&=&(q-q^{-1})\sum_{s>i}x^k_sx^s_j 
\qquad\mbox{ for }j<i<k;\\  
x^i_jx^k_l-x^k_lx^i_j&=&(q-q^{-1})x^k_jx^i_l 
\qquad\mbox{ for }i>k,j<l,i\neq l, j\neq k, i\neq j;\\  
x^i_jx^k_l-x^k_lx^i_j&=&(q-q^{-1})x^i_lx^k_j 
\qquad\mbox{ for }i>k,j<l,i\neq l,j\neq k,k\neq l;\\ 
x^j_lx^i_j-qx^i_jx^j_l&=&(q-q^{-1})(-x^j_jx^i_l+\sum_{s>j}x^i_sx^s_l) 
\qquad\mbox{ for }j<l, j<i, i\neq l;\\  
x^i_jx^k_i-qx^k_ix^i_j&=&(q-q^{-1})(x^k_jx^i_i+\sum_{s>i}x^k_sx^s_j) 
\qquad\mbox{ for }k<i, j<i, k\neq j;\\  
x^i_jx^j_i-x^j_ix^i_j
&=&(1-q^{-2})(x^j_jx^i_i-\sum_{s>j}x^i_sx^s_i
+\sum_{t>i}x^j_tx^t_j) 
\qquad\mbox{ for }j<i. 
\end{eqnarray*} 

\begin{proposition}\label{prop:noether} 
\begin{itemize} 
\item[(i)] The algebra $\rea$ is noetherian. 
\item[(ii)] Let $x_1,\ldots,x_{N^2}$ be an arbitrary permutation of 
$x^1_1,\ldots,x^N_N$. 
Then the monomials of the form 
$x^{a_1}_1\cdots x^{a_{N^2}}_{N^2}$ constitute a $\mc$-vector space basis 
of $\rea$. 
\item[(iii)] $\rea$ is a domain. 
\end{itemize} 
\end{proposition} 

\begin{proof} (i) Denote by $I^i_j$ the ideal of $\rea$ generated by 
the $x^s_t$ with $s\geq i$, $t\geq j$, $(s,t)\neq (i,j)$. 
The relation indexed by $i,j,k,l$ in \eqref{eq:rea-coord} says that 
for some $q$-power $q^{ik}_{jl}$, we have 
$x^i_jx^k_l-q^{ik}_{jl}x^k_lx^i_j\in I^i_j\cap I^k_l$. 
In particular, since $I^N_N=0$, the element $x^N_N$ is normal in $\rea$; 
that is, $x^N_N\rea=\rea x^N_N$. 
Next, observe that $I^{N-1}_N=\langle x^N_N\rangle=I^N_{N-1}$. 
It follows that both the images of $x^N_{N-1}$ and  $x^{N-1}_N$ are 
normal in the quotient $\rea\langle x^N_N\rangle$. 
Continuing this way it is easy to see that the generators can be arranged 
into a so-called polynormal sequence. 
To be more precise, let $x_1,\ldots,x_{N^2}$ be a permutation of 
$x^1_1,\ldots,x^N_N$. Define the function $p$ by $x_{p(i,j)}=x^i_j$. 
Assume that 
\begin{equation}\label{eq:cond-polynorm} 
\mbox{for all }i,j,k,l \mbox{ with }i\leq k, j\leq l 
\mbox{ we have that }p(k,l)\leq p(i,j).
\end{equation}  
Then $x_1,\ldots,x_{N^2}$ is a polynormal sequence: 
for all $i$, the image of $x_i$ is normal in 
$\rea/\langle x_1,\ldots,x_{i-1}\rangle$.  
Condition \eqref{eq:cond-polynorm} holds for example for the sequence 
$x^N_N,x^{N-1}_N,\ldots,x^1_N,x^N_{N-1},x^{N-1}_{N-1},\ldots,x^1_1$. 
So the graded algebra $\rea$ is generated by a polynormal sequence, hence is 
noetherian (both left and right) by Lemma 8.2 of \cite{atvdb}. 
(Alternatively, one may apply Proposition I.8.17 from \cite{bg} for a 
generating sequence satisfying \eqref{eq:cond-polynorm} to conclude 
that $\rea$ is noetherian.) 

(ii) Write $W$ for the free semigroup generated by the variables $x^i_j$. 
Define the function $h$ on $W$ by 
\[h(x^i_j\cdots x^k_l)=ij+\cdots+kl.\]  
Note that if $w'$ is obtained by permuting the factors of the monomial $w$, 
then $h(w')=h(w)$. 
Introduce a partial order $\prec$ on $W$: for $v,w\in W$, we set  
$v\prec w$ if $h(v)>h(w)$. Clearly, $\prec$ is compatible with multiplication 
in the sense that $v\prec w$ implies $avb\prec awb$ for all 
$a,b\in W\cup\{1\}$. 
In $\rea$, the relation indexed by $i,j,k,l$ in \eqref{eq:rea-coord} 
says that 
$x^i_jx^k_l-q^{ik}_{jl}x^k_lx^i_j$ is a linear combination of monomials $w$ 
with $w\prec x^i_jx^k_l$ and $w\prec x^k_lx^i_j$ 
(where $q^{ik}_{jl}$ is a $q$-power, hence non-zero). 
Restrict attention to the degree $d$ homogeneous component $\rea^{(d)}$ 
of $\rea$. By an iterated use of the quadratic relations, 
we can rewrite 
$w=x^{i_1}_{j_1}\cdots x^{i_d}_{j_d}$ as a non-zero scalar multiple of 
some $w'=(x_1)^{a_1}\cdots (x_{N^2})^{a_{N^2}}$ 
(obtained by permuting the factors of $w$), 
modulo the span of monomials $v$ with $v\prec w$, $v\prec w'$. 
This obviously implies that the monomials of the form 
$(x_1)^{a_1}\cdots (x_{N^2})^{a_{N^2}}$ span $\rea^{(d)}$. 
Their number coincides with the dimension of $\rea^{(d)}$ 
(which we know already), so they form a basis. 

(iii) Given a homogeneous element $f$ of $\rea$, 
express it in terms of the basis described in (ii). 
Among the monomials of $f$ that are maximal with respect to the partial 
ordering $\prec$ choose the one which is largest with respect to the
lexicographic ordering of the monomials of the form 
$(x_1)^{i_1}\cdots (x_{N^2})^{i_{N^2}}$, and call it the 
{\it leading term of} $f$. 
If $f,g$ are non-zero homogeneous elements of $\rea$ with leading term 
$(x_1)^{i_1}\cdots (x_{N^2})^{i_{N^2}}$ and 
$(x_1)^{j_1}\cdots (x_{N^2})^{j_{N^2}}$,  
then by the rewriting process sketched in the proof of (ii), the leading term 
of $fg$ is 
$(x_1^{i_1+j_1})\cdots (x_{N^2})^{i_{N^2}+j_{N^2}}$, so 
$fg$ is non-zero. Thus $\rea$ does not contain zero-divisors. 
\end{proof} 


\section{Quantum trace rings}\label{sec:q-trace-ring} 

After these preparations we are in position to define 
$\qmatinv$ as the subalgebra of $\mcoact$-coinvariants 
in $\qsevmat$; that is, 
\[\qmatinv=\{f\in\qsevmat\mid \mcoact f=f\otimes 1\}.\] 
Let $\qgenmat$ denote the unitary subalgebra of $M_N(\qsevmat)$ 
(the algebra of $N\times N$ matrices over $\qsevmat$) 
generated 
by $X(1),\ldots,X(m)$. And finally, consider 
\[\qtracering
=\{F\in M_N(\qsevmat)\mid \mcoact F=S(T)FT\},\] 
where $\mcoact F$ stands for the matrix obtained by applying $\mcoact$ 
to each entry of $F$, 
and $S(T)XT$ is identified with an element of 
$M_N(\qsevmat)\otimes\qsl$ in the way indicated in 
Section~\ref{sec:quadratic}. 

\begin{proposition}\label{prop:subalgebra}
$\qtracering$ is a subalgebra of $M_N(\qsevmat)$. 
\end{proposition} 

\begin{proof} 
This follows from the formula $S(T)=T^{-1}$, together with 
the multiplicativity of the coaction 
$\mcoact:\qsevmat\to\qsevmat\otimes\qsl$.  
\end{proof} 

The generic matrices $X(k)$ belong to $\qtracering$ by definition of 
$\mcoact$, hence all products $X(k_1)\cdots X(k_d)$ belong to 
$\qtracering$ by Proposition~\ref{prop:subalgebra}. 
Also, if $f\in\qmatinv$, then the scalar matrix $fI$ is contained 
in $\qtracering$. 

One can produce $\mcoact$-coinvariants using the 
so-called {\it quantum trace} operation.  
For an arbitrary $N\times N$ matrix $F$ with entries from a 
$\mc$-algebra, define 
$\qtr(F)=\Tr(QF)$, 
where $Q$ is the diagonal matrix with diagonal entries 
$q^{N-1},q^{N-3},q^{N-5},\ldots,q^{-N+1}$. 
In particular, $\qtr(I)=[N]_q$, the $N$th $q$-integer. 
Now $\qtr(X)$ is a $\coact$-coinvariant in $\rea$; more generally,   
if $F\in\qtracering$, then $\qtr(F)$ is contained in $\qmatinv$ 
(see the proof of Theorem~\ref{th:1}). 
In particular, 
the elements $\qtr(X(k_1)\cdots X(k_d))$ are all contained 
in $\qmatinv$.  

Note that there is a natural $\mz^m$-grading on the algebra $\qsevmat$: 
the multidegree of the entries of $X(k)$ is the $k$th standard basis 
vector of $\mz^m$. This induces a $\mz^m$-grading on $M_N(\qsevmat)$ as well. 
Since the multihomogeneous components of $\qsevmat$ are subcomodules 
with respect to $\mcoact$, 
we have that $\qmatinv$ is a $\mz^m$-graded subalgebra of $\qsevmat$, 
and similarly, $\qtracering$ is a $\mz^m$-graded subalgebra of 
$M_N(\qsevmat)$. 

\begin{theorem}\label{th:1}
Assume that $q\in\mc^{\times}$ is not a root of unity. Then 
the $\mz^m$-graded algebra $\qmatinv$ has the same Hilbert series 
as its classical counterpart $\matinv$, 
and $\qtracering$ has the same Hilbert series as 
$\tracering$. 
\end{theorem} 

\begin{proof} 
Consider the multihomogeneous component 
$\qsevmat^{(d_1,\ldots,d_m)}$, it is a subcomodule of $\qsevmat$ 
with respect to $\coact$. Since $q$ is assumed to be not a root of unity, 
$\qsl$ is cosemisimple (as in the case that $q=1$), 
and this finite dimensional comodule decomposes as the 
direct sum of simple subcomodules. By definition, the coefficient of 
$t_1^{d_1}\cdots t_m^{d_m}$ in the Hilbert series 
$H(\qmatinv;t_1,\ldots,t_m)$ is the multiplicity of the trivial comodule as a 
direct summand. 
This multiplicity can be computed by 
``restricting the coaction to the diagonal subgroup''  
(see for example \cite{dl} for a more detailed explanation). 
To be more precise, 
denote by $\fundiag$ the coordinate Hopf algebra of the diagonal subgroup 
of $SL_N(\mc)$, so 
$\fundiag=\mc[z_1,\ldots,z_N\mid z_1\cdots z_N=1]$. 
There is a surjective homomorphism $\pi:\qsl\to\fundiag$, 
$\pi(t^i_i)=z_i$, and $\pi(t^i_j)=0$ for $i\neq j$. 
Consider the coaction 
$(\id\otimes\pi)\mcoact:\qsevmat\to\qsevmat\otimes\fundiag$ 
of $\fundiag$. It follows from Proposition~\ref{prop:noether} (ii) 
and the interpretation of $\qsevmat$ as a braided tensor product of copies of 
$\rea$ that 
the $\fundiag$-comodule $\qsevmat^{(d_1,\ldots,d_m)}$ 
is isomorphic to its counterpart in the classical case $q=1$. 
That is, formally the same computation gives the coefficients of the 
Hilbert series of $\qmatinv$ as that of $\matinv$. 

The proof of the statement about the Hilbert series of $\qtracering$ 
is similar, one just needs an appropriate interpretation of 
the elements of $\qtracering$ 
in terms of corepresentation theory. 
Denote by $\varphi$ the restriction of the coaction $\beta$ 
to the $N^2$-dimensional space $\xijspan=\linspan\{x^i_j\}_{i,j=1}^N$. 
As in the classical case $q=1$, 
this corepresentation decomposes as $\varphi\cong \triv\oplus \varphi_0$, 
where $\triv$ is the trivial corepresentation, and 
$\varphi_0$ is an irreducible $(N^2-1)$-dimensional corepresentation. 
It is well known (and easy to show; we shall comment on it later in this 
Section) that $\triv$ is realized on the 
subspace of $\xijspan$ spanned by $\qtr(X)$.  
Thus the decomposition of $\xijspan$ as a sum of simple subcomodules 
is $\xijspan=\mc\qtr(X)\oplus\xijspan_0$, 
where $\xijspan_0$ is 
spanned by the entries of $X-\frac{\qtr(X)}{[N]_q}I$, 
with $[N]_q$ standing for the $q$-integer 
$q^{N-1}+q^{N-3}+\cdots+q^{-N+1}$. 

By definition, a matrix $F=(f^i_j)$ belongs to  
$\qtracering$ if and only if its  
entries span a $\mcoact$-subcomodule of $\qsevmat$, 
such that $\xijspan\to\linspan\{f^i_j\}$, $x^i_j\mapsto f^i_j$  
is a homomorphism of $\qsl$-comodules. 
It follows in particular that if $F\in\qtracering$, 
then $\qtr(F)\in\qmatinv$. 
Therefore $G=F-\frac{\qtr(F)}{[N_q]}I\in\qtracering$ with 
$\qtr(G)=0$, and the entries of $G$ span a subcomodule of 
$\qsevmat$ which is a homomorphic image of the simple comodule $\xijspan_0$. 
Conversely, given an element $f\in\qmatinv$, the scalar matrix 
$fI$ obviously belongs to $\qtracering$. 
Now take a subcomodule $\mathcal{Y}$ 
of $\qsevmat$ which is isomorphic to $\xijspan_0$. 
Fix a comodule surjection $\xijspan\to\mathcal{Y}$, $x^i_j\mapsto f^i_j$. 
(Note that by Schur's Lemma, this homomorphism is determined up to a scalar 
multiple.) 
Then the matrix $F=(f^i_j)$ belongs to $\qtracering$ by construction, 
and has the property $\qtr(F)=0$. 

Putting this all together, we reach the conclusion that 
the dimension of $\qtracering^{(d_1,\ldots,d_m)}$ is the sum of 
the dimension of $\qmatinv^{(d_1,\ldots,d_m)}$ and the multiplicity 
of the irreducible corepresentation $\varphi_0$ as a direct summand 
of $\mcoact$ on 
$\qsevmat^{(d_1,\ldots,d_m)}$. This latter multiplicity can be computed 
by restriction to the corresponding $\fundiag$-coaction, 
which by Proposition~\ref{prop:noether} (ii) is isomorphic to its classical 
counterpart. Finally, note that the dimension of 
$\tracering^{(d_1,\ldots,d_m)}$ in the classical case 
clearly can be expressed in the 
same way in terms of multiplicities of two types of 
irreducible summands in $\sevmat$. 
\end{proof}

\begin{remark}\label{rem:1}
{\rm Let us extract from the above proof explicitly that the $\mc$-vector
space $\qtracering$ can be identified with the space of 
$\qsl$-comodule homomorphisms 
$\Homom_{\qsl}(\xijspan,\qsevmat)$. 
There is a natural vector space isomorphism 
$\Homom_{\mc}(\xijspan,\qsevmat)\cong\qsevmat\otimes\xijspan^*$. 
The dual space $\xijspan^*$ of $\xijspan$ is an $\qsl$-comodule via the
contragredient corepresentation of $\varphi$ (cf. section 11.1.3, page 398 
in \cite{kl-sch}), and consider the tensor product comodule structure on 
$\qsevmat\otimes\xijspan^*$. 
It is straightforward to check that under the isomorphism 
$\Homom_{\mc}(\xijspan,\qsevmat)\cong\qsevmat\otimes\xijspan^*$ 
the space of coinvariants $(\qsevmat\otimes\xijspan^*)^{\qsl}$ is 
mapped onto $\Homom_{\qsl}(\xijspan,\qsevmat)$. 
In particular, we have 
$\qtracering\cong (\qsevmat\otimes\xijspan^*)^{\qsl}$ 
as vector spaces. } 
\end{remark} 

Under a stronger restriction on $q$ we are able to present explicit 
generators of $\qmatinv$ and $\qtracering$. 

\begin{theorem}\label{th:2} 
Assume that $q$ is transcendental over the rationals. Then 
$\qmatinv$ is the $\mc$-subalgebra of $\qsevmat$ 
generated by the elements of the form 
$\qtr(X(i_1)\cdots X(i_d))$, 
where $i_1,\ldots,i_d\in\{1,\ldots,m\}$, and $d\leq N^2$. 
Furthermore, $\qtracering$ is the left (respectively, right) 
$\qmatinv$-submodule of 
$M_N(\qsevmat)$ generated by the elements 
$X(i_1)\cdots X(i_d)$ with $d\leq N^2-1$. 
\end{theorem} 

\begin{proof} 
These statements are known to hold when $q=1$. 
Take a set of multihomogeneous 
monomials whose traces form a basis of $\matinv$ when $q=1$. 
Take the $q$-traces of formally the same monomials (but interpret them  
as matrices with entries from $\qsevmat$). We obtain some multihomogeneous 
elements in $\qmatinv$. Since $q$ is assumed to be transcendental, 
we may conclude by a standard argument that 
these elements in $\qmatinv$ are linearly independent
(because of the linear independency of the corresponding elements for $q=1$). 
Therefore they constitute a basis of $\qmatinv$ by Theorem~\ref{th:1}. 
The proof of the statement about $\qtracering$ is the same. 
\end{proof} 

When $q$ is not a root of unity, 
the multilinear component $\qsevmat^{(1,\ldots,1)}$ 
can be naturally identified with a homomorphic image of 
the Hecke algebra of the symmetric group $S_m$. 
This follows from some natural linear algebra 
isomorphisms (mentioned already in Remark~\ref{rem:1}). 
Indeed, consider the fundamental corepresentation 
$\omega:V\to V\otimes\qsl$, $e_j\mapsto e_i\otimes t^i_j$, 
where $V$ is an $N$-dimensional vector space with basis $e_1,\ldots,e_N$. 
Write $\varepsilon^1,\ldots,\varepsilon^N$ for the corresponding 
dual basis in $V^*$. 
Then using the notation of the proof of Theorem~\ref{th:1}, 
$\xijspan\to V^*\otimes V$, $x^i_j\mapsto \varepsilon^i\otimes e_j$ 
intertwines $\varphi$ and $\omega^*\otimes \omega$. 
Therefore as a comodule, the multilinear component 
$\qsevmat^{(1,\ldots,1)}$ can be identified with 
$V^*\otimes V\otimes\cdots \otimes V^*\otimes V$. 
With a repeated use of the comodule isomorphism 
$V\otimes V^*\to V^*\otimes V$, 
$e_i\otimes\varepsilon^j\mapsto \rtilde^{kj}_{il}\varepsilon^l\otimes e_k$, 
and the other natural isomorphisms mentioned in 
Remark~\ref{rem:1} (note also that 
$(A\otimes B)^*\cong B^*\otimes A^*$ for any finite dimensional 
$\qsl$-comodules $A,B$) 
we obtain an isomorphism 
$\qsevmat^{(1,\ldots,1)}\cong \Endom_{\mc}(V\otimes\cdots\otimes V)$, 
which restricts to 
$\qmatinv^{(1,\ldots,1)}\cong\Endom_{\qsl}(V^{\otimes m})$. 
The latter is known to be generated as an algebra by $\rhat_{i,i+1}$ 
(operating as $\rhat$ on the $i$th and $(i+1)$th tensor component, and 
operating as the identity on the other components), $i=1,\ldots,m-1$. 
That way we get an explicit vector space homomorphism from the Hecke algebra 
of the symmetric group $S_m$ onto $\qmatinv^{(1,\ldots,1)}$. 
The kernel of this homomorphism is a known two-sided ideal of the 
Hecke algebra. 
For example, in the special case $m=1$, the element 
$\id\in\Endom_{\qsl}(V)$ is mapped to $q^{-2}\qtr(X)$. 
However, it does not seem to be easy to exploit this 
connection for explicit calculations, because of the excessive use of the 
R-matrix in some of the above identifications. 
For example, in the case $N=m=2$, we have 
\[\id_{V\otimes V}\mapsto 
e_i\otimes e_j\otimes\varepsilon^j\otimes\varepsilon^i\mapsto 
\rtilde^{si}_{kt}\rtilde^{ml}_{in}\rtilde^{kj}_{jl}x^n_my^s_t,\] 
which is 
$q^{-2}\qtr(XY)+(q^{-5}-q^{-3})\qtr(X)\qtr(Y)\in\twotwoqsevmat$. 

Another interesting fact about $\qsevmat$ is the centrality of certain 
$\mcoact$-coinvariants. 
The $\coact$-coinvariants are central in $\rea$ (see e.g. \cite{majid2}),  
and $\qtr(X(k)^n)$ is central in $\qsevmat$, as it was 
observed in \cite{ku-sa}. Using the approach of \cite{majid2} it is easy 
to  show that 
if $f\in\qmatinv$ belongs to the subalgebra generated by the entries 
of $X(k)$ for some $k$, then $f$ is central in $\qsevmat$. 
It is sufficient to show that $f$ commutes with the entries of $X(l)$ 
for $l\neq k$ (since we know already that $f$ commutes with the entries of 
$X(k)$). So one can reduce the question to the case $m=2$. 
Let $x,y$ be elements of $\ntwoqsevmat$, $x$ depending only on the entries 
of $X$, and $y$ depending only on the entries of $Y$. 
With Sweedler's notation, $\coact(x)=\sum x_0\otimes x_1$, 
and $\coact(y)=\sum y_0\otimes y_1$. 
By definition of the multiplication in the braided tensor product 
$\rea\btensor\rea$, we have 
\begin{equation}\label{eq:braided-tensor}
yx=\sum r(y_1,x_1)x_0y_0,
\end{equation} 
where $r$ is the universal r-form 
in the coquasitriangular Hopf algebra $\qsl$ 
(see e.g. Proposition 35, 10.3.2 in \cite{kl-sch}). 
If $\coact(x)=x\otimes 1$, then we get 
$yx=\sum r(y_1,1)xy_0=x\sum \varepsilon(y_1)y_0=xy$. 
Similarly, $\coact(y)=y\otimes 1$ implies $yx=xy$. 

We mention also that the entries of $XY$ satisfy the defining relations of 
$\rea$; 
this is stated in \cite{ku-sa}. 
Indeed, 
\begin{eqnarray*} 
(XY)_2R_{21}(XY)_1R &=&X_2Y_2\tau R\tau X_1Y_1R=
X_2\tau (Y_1RX_2)\tau Y_1R \\ 
&=&X_2\tau RX_2R^{-1}Y_1R\tau Y_1 R
=(X_2R_{21}X_1)\tau R^{-1}\tau (Y_2R_{21}Y_1R)
\\ 
&=&R_{21}X_1RX_2R^{-1}R_{21}^{-1}R_{21}Y_1RY_2
=R_{21}X_1R(X_2R^{-1}Y_1R)Y_2
\\ &=&R_{21}X_1RR^{-1}Y_1RX_2Y_2
=R_{21}X_1Y_1RX_2Y_2
\\ &=& R_{21}{(XY)}_1R{(XY)}_2, 
\end{eqnarray*}
so $XY$ satisfies the matrix equality \eqref{eq:rea}. 
Furthermore, the map $X\mapsto X(1)X(2)$, $Y\mapsto X(3)$ 
extends to an algebra homomorphism $\ntwoqsevmat\to\nthreeqsevmat$. 
All one needs to check is that the matrices $X(1)X(2)$ and $X(3)$ 
satisfy the equality \eqref{eq:qsevmatntwo}. 
This can be done as in the calculation above. 
Alternatively, think of $\nthreeqsevmat$ as $\ntwoqsevmat\btensor\rea$; 
the multiplication in the braided tensor product is expressed in terms 
of the $\qsl$-comodule algebra structures of the factors 
(see \eqref{eq:braided-tensor}). Since the linear isomorphisms 
$X\mapsto X(1)X(2)$, $Y\mapsto X(3)$ intertwine the corresponding coactions, 
the cross relations 
\eqref{eq:qsevmatntwo-coord} hold between the entries of $X(1)X(2)$ and $X(3)$.


\section{The $2\times 2$ case}\label{sec:2x2} 

Throughout this section we assume $N=2$. 
Recall from Section~\ref{seq:polynormal} 
that $\tworea$ is generated by $x^1_1,x^1_2,x^2_1,x^2_2$, subject to 
the relations 
\begin{align} \label{eq:2rea} 
&x^2_2x^1_2=q^2x^1_2x^2_2;&\qquad
&x^1_1x^1_2=x^1_2x^1_1+(q^{-2}-1)x^1_2x^2_2;&\\ \notag
&x^1_1x^2_2=x^2_2x^1_1;&\qquad
&x^2_1x^1_2=x^1_2x^2_1+(q^{-2}-1)x^2_2(x^2_2-x^1_1);&\\ \notag
&x^2_1x^2_2=q^2x^2_2x^2_1;&\qquad
&x^2_1x^1_1=x^1_1x^2_1+(q^{-2}-1)x^2_2x^2_1.&
\end{align} 

The cross relations in $\twotwoqsevmat$ between the entries of $X$ and $Y$ are 
the following: 

\begin{eqnarray}\label{eq:cross} 
y^1_1x^1_1&=&x^1_1y^1_1+(q^{-4}-q^{-2})x^2_1y^1_2;\\ \notag
y^1_2x^1_1&=&x^1_1y^1_2;\\ \notag 
y^1_1x^1_2&=&x^1_2y^1_1+(1-q^{-2})x^1_1y^1_2+(q^{-2}-1)x^2_2y^1_2;\\ \notag
y^1_2x^1_2&=&q^2x^1_2y^1_2;\\ \notag
y^2_1x^1_1&=&x^1_1y^2_1+(q^{-2}-1)x^2_1y^2_2+(1-q^{-2})x^2_1y^1_1;\\ \notag
y^2_2x^1_1&=&x^1_1y^2_2+(1-q^{-2})x^2_1y^1_2;\\ \notag
y^2_1x^1_2&=&q^{-2}x^1_2y^2_1+(q^{-2}-1)x^1_1y^1_1+(1-q^{-2})x^1_1y^2_2
+(1-q^{-2})x^2_2y^1_1+
\\ \notag &~&+(q^{-2}-1)x^2_2y^2_2+(q^{-4}-q^{-2}-1+q^2)x^2_1y^1_2;\\ \notag
y^2_2x^1_2&=&x^1_2y^2_2+(1-q^2)x^1_1y^1_2+(q^2-1)x^2_2y^1_2;\\ \notag
y^1_1x^2_1&=&x^2_1y^1_1;\\ \notag
y^1_2x^2_1&=&q^{-2}x^2_1y^1_2;\\ \notag
y^1_1x^2_2&=&x^2_2y^1_1+(1-q^{-2})x^2_1y^1_2;\\ \notag
y^1_2x^2_2&=&x^2_2y^1_2;\\ \notag
y^2_1x^2_1&=&q^2x^2_1y^2_1;\\ \notag
y^2_2x^2_1&=&x^2_1y^2_2;\\ \notag
y^2_1x^2_2&=&x^2_2y^2_1+(q^2-1)x^2_1y^2_2+(1-q^2)x^2_1y^1_1;\\ \notag
y^2_2x^2_2&=&x^2_2y^2_2+(1-q^2)x^2_1y^1_2.
\end{eqnarray} 

The algebra $\twomqsevmat$ is generated by the entries of the $2\times 2$ 
matrices $X(1),\ldots,X(m)$. For each $k$, the entries of $X(k)$ satisfy 
\eqref{eq:2rea} with $X\mapsto X(k)$, and 
for each $k<l$, the entries of $X(k)$ and $X(l)$ satisfy 
\eqref{eq:cross} with $X\mapsto X(k)$, $Y\mapsto X(l)$. 
The monomials  
\begin{equation} \label{eq:qsevmat22-basis}
w(1)\ldots w(m)\quad\mbox{ with }\quad 
w(k)=(x(k)^1_2)^{a_k}(x(k)^2_2)^{b_k}(x(k)^1_1)^{c_k}(x(k)^2_1)^{d_k}
\end{equation} 
form a basis of $\twomqsevmat$. Above we have presented the relations of 
$\twomqsevmat$ in such a form that makes obvious an algorithm to 
rewrite an arbitrary monomial in the generators as a linear combination of 
these basis elements. This allows us to perform calculations in $\twomqsevmat$. 

The quantum version of the Cayley-Hamilton identity for the matrix 
$X$ of the generators of $\rea$ 
was found in \cite{ps} (see also \cite{iop}). 
In the special case $N=2$ it takes the form 
\begin{equation}\label{eq:CH}
X^2-q^{-1}\qtr(X)X
+\frac{1}{[2]_q}(q^{-1}\qtr(X)^2-\qtr(X^2))I=0.
\end{equation} 
We note that 
$\frac{1}{q+q^{-1}}(q^{-1}\qtr(X)^2-\qtr(X^2))
=q^{-2}x^2_2x^1_1-x^1_2x^2_1$. 

Classically one can polarize the $2\times 2$ 
Cayley-Hamilton identity to obtain an equivalent 
bilinear identity depending on two matrix variables. 
This bilinear identity has an analogues in the quantum setting: 

\begin{proposition}\label{prop:bilin-CH} 
The following equality holds in $\twotwoqtracering$: 
\begin{eqnarray}\label{eq:bilin-CH} 
XY+q^2YX-q\qtr(X)Y-q\qtr(Y)X
\notag \\ +(\qtr(X)\qtr(Y)-q^{-1}\qtr(XY))I=0.
\end{eqnarray}  
\end{proposition} 

\begin{proof} 
Direct computation, using the above mentioned basis and rewriting 
algorithm in $\twotwoqsevmat$. 
\end{proof} 

We take a digression and for illustrative purposes we list a couple of 
other equalities in $\twotwoqsevmat$ and in $\twothreeqsevmat$. 
We have 
\begin{equation}\label{eq:qtryx}
q^2\qtr(YX)-q^{-2}\qtr(XY)=(q-q^{-1})\qtr(X)\qtr(Y) 
\mbox{ in }\twotwoqsevmat
\end{equation}  
(this can be checked by direct computation again). 
To simplify notation, write $X,Y,Z$ (instead of $X(1),X(2),X(3)$) 
for the matrices of the generators of 
$\twothreeqsevmat$. 
Multiplying \eqref{eq:bilin-CH} by $Z$ from the right and taking $q$-trace 
one obtains 
\begin{eqnarray}\label{eq:fundqtrid2} 
\qtr(X)\qtr(Y)\qtr(Z)-q^{-1}\qtr(XY)\qtr(Z)
-q\qtr(X)\qtr(YZ) \notag
\\ -q\qtr(Y)\qtr(XZ)+q^2\qtr(YXZ)+\qtr(XYZ)=0 
\mbox{ in }\twothreeqsevmat,   
\end{eqnarray} 
whereas making the substitution $X\mapsto Y$, $Y\mapsto Z$ in 
\eqref{eq:bilin-CH}, multiplying by $X$ from the left, and taking 
$q$-trace one gets 
\begin{eqnarray}\label{eq:fundqtrid1}
\qtr(X)\qtr(Y)\qtr(Z)-q\qtr(XY)\qtr(Z)
-q^{-1}\qtr(X)\qtr(YZ) \notag
\\ -q\qtr(XZ)\qtr(Y)
+q^2\qtr(XZY)+\qtr(XYZ)=0 
\mbox{ in }\twothreeqsevmat.  
\end{eqnarray}
From either of the above two equalities one recovers the so-called 
fundamental trace 
identity of $2\times 2$ matrices as the special case $q=1$. 
However, we have no canonical choice to single out one of them as 
``the fundamental $q$-trace identity'', because as is indicated by 
\eqref{eq:qtryx}, 
the number of the different trilinear products of $q$-traces 
is larger than $6$, whereas $\dim_{\mc}(\twothreeqsevmat^{(1,1,1)})=5$, 
as in the classical case. 
Furthermore, make the substitutions $X\mapsto XY$, $Y\mapsto Z$ 
(respectively $X\mapsto X$, $Y\mapsto YZ$) in the 
identity \eqref{eq:qtryx}; 
by the discussion at the end of Section~\ref{sec:q-trace-ring}, we obtain 
the equalities  
\begin{equation}\label{eq:qtrzxy}
q^2\qtr(ZXY)-q^{-2}\qtr(XYZ)=(q-q^{-1})\qtr(XY)\qtr(Z),  
\end{equation} 
\begin{equation}\label{eq:qtryzx}
q^2\qtr(YZX)-q^{-2}\qtr(XYZ)=(q-q^{-1})\qtr(X)\qtr(YZ).   
\end{equation} 
Finally, multiply \eqref{eq:bilin-CH} by $Z$ from the left and take its 
$q$-trace; from this, \eqref{eq:fundqtrid1}, \eqref{eq:fundqtrid2}, 
\eqref{eq:qtrzxy}, and \eqref{eq:qtryzx}
one gets that 
modulo the subalgebra generated by q-traces of monomials of degree $\leq 2$, 
we have 
$\qtr(XYZ)\equiv-q^2\qtr(YXZ)\equiv -q^2\qtr(XZY)
\equiv q^4\qtr(YZX)\equiv q^4\qtr(ZXY)\equiv -q^6\qtr(ZYX)$. 

\bigskip
Next we investigate $\twotwoqtracering$. 
In order to simplify notation, set 
\begin{eqnarray*}
A=q^{-1}\qtr(X)I;\quad B=q^{-1}\qtr(Y)I;
\\ C=(q^{-2}x^2_2x^1_1-x^1_2x^2_1)I; \quad
D=(q^{-2}y^2_2y^1_1-y^1_2y^2_1)I;
\\E=AB-q^{-3}\qtr(XY)I.
\end{eqnarray*} 

\begin{theorem}\label{th:presentation} 
Suppose that $q$ is transcendental over the rationals. 
Then we have the following. 
\begin{itemize}
\item[(i)] $\twotwoqmatinv$ is a $5$-variable commutative polynomial algebra 
generated by $\qtr(X),\qtr(X^2),\qtr(Y),\qtr(Y^2),\qtr(XY)$. 
\item[(ii)] $\twotwoqtracering$ is a free left $\twotwoqmatinv$-module 
generated by $I,X,Y,XY$. 
\item[(iii)] As a $\mc$-algebra, $\twotwoqtracering$ is generated by 
$X,Y,A,B,C,D,E$, and a complete list of relations among these generators is 
the following: 
\begin{eqnarray*} 
A,B,C,D \mbox{ are central};
\\ X^2=AX-C;\quad Y^2=BY-D;
\\ YX=-q^{-2}XY+AY+BX-E;
\\ XE=q^2EX+(1-q^2)ABX+(q^{-2}-q^2)CY+
\\ +(1-q^{-2})AXY+(q^2-1)CB;
\\ YE=q^{-2}EY+(1-q^{-2})ABY+(1-q^{-4})DX+
\\ +(q^{-4}-q^{-2})BXY+(q^{-2}-1)AD.
\end{eqnarray*} 
\end{itemize}
\end{theorem}

\begin{proof} 
(i) As we noted in Section~\ref{sec:q-trace-ring}, 
$\qtr(X),\qtr(X^2),\qtr(Y),\qtr(Y^2)$ are central in $\twotwoqsevmat$, 
hence the given five elements pairwise commute. They are algebraically 
independent because of the algebraic independence of the corresponding 
elements in the case $q=1$ (here we use the assumption on $q$, in the same 
way as in the proof of Theorem~\ref{th:2}). 
We know from 
Theorem~\ref{th:1} that $\twotwoqmatinv$ has the same Hilbert series 
as a polynomial algebra generated by two degree $1$ elements and three 
degree $2$ elements. 
This implies that $\twotwoqmatinv$ coincides with its subalgebra 
$\mc[\qtr(X),\qtr(X^2),\qtr(Y),\qtr(Y^2),\qtr(XY)]$.  

(ii) $\twotwoqtracering$ contains $I,X,Y,XY$, hence the left $\twotwoqmatinv$-module 
generated by them. This is a free module, because the corresponding 
module is free when $q=1$ (again the transcendentality of $q$ is used). 
Consequently, it has the same Hilbert series as $\twotwotracering$ 
(computed in \cite{fhl}). Hence by Theorem~\ref{th:1}, 
the left $\twotwoqmatinv$-module generated by $I,X,Y,XY$ 
coincides with $\twotwoqtracering$. 

(iii) It follows from (ii) that $X,Y,A,B,C,D,E$ generate $\twotwoqtracering$ 
as an algebra. By the remarks in Section~\ref{sec:q-trace-ring} 
about the centrality of certain $\mcoact$-coinvariants, we have that $A,B,C,D$ 
are central in $\twotwoqtracering$. 
The relations $X^2=AX-C$, $Y^2=BY-D$ are the Cayley-Hamilton identities 
of \cite{ps}, see \eqref{eq:CH}. The next relation in the list 
is the bilinear Cayley-Hamilton identity 
of Proposition~\ref{prop:bilin-CH}. 
The remaining two relations are verified by direct calculation, using the 
basis and rewriting algorithm in $\twotwoqsevmat$ explained at the beginning 
of this Section. 
Thus the relations we listed all hold. It is straightforward 
that modulo these relations an arbitrary product of the elements 
$X,Y,A,B,C,D,E$ can be rewritten as an element of the left 
$\mc[A,B,C,D,E]$-module generated by $I,X,Y,XY$. 
Hence by (ii), the given list of relations is complete. 
\end{proof} 

\begin{remark} 
{\rm A rather lengthy calculation using the relations in (iii) 
yields that $XYE=EXY$. Knowing that $A,B$ are central, this is equivalent 
to the assertion that $\qtr(XY)$ commutes with $XY$. 
However, this latter statement immediately follows (for arbitrary $N$) 
from the existence of the algebra homomorphism $\rea\to\ntwoqsevmat$, 
$X\mapsto XY$ (discussed in 
Section~\ref{sec:q-trace-ring}), and the centrality of 
$\qtr(X)$ in $\rea$. }
\end{remark} 

Now consider the ordinary trace ring $\twotwotracering$, 
and write $x,y$ for the matrix generators, so $x,y$ are 
$2\times 2$ generic matrices with pairwise 
commuting algebraically independent entries. 
The algebra $\twotwotracering$ probably appeared first in \cite{syl};  
see \cite{f2} for a recent survey.  
We set  
\begin{eqnarray*} a=\Tr(x)I;\quad b=\Tr(y)I;\quad 
c=\det(x)I; \quad d=\det(y)I;
\\ e=ab-\Tr(xy)I.
\end{eqnarray*} 

\begin{theorem}\label{th:isomorphism} 
Suppose that $q$ is transcendental over the rationals. 
Then the map 
\begin{eqnarray*} X\mapsto x,\ Y\mapsto y,\ A\mapsto a, B\mapsto b,\ 
C\mapsto c,\ D\mapsto d,
\\ E\mapsto e+(1-q^{-2})xy
\end{eqnarray*} 
extends to an algebra isomorphism 
$\twotwoqtracering\to\twotwotracering$. 
In particular, we have 
$\twotwoqgenmat\cong\twotwogenmat$ are isomorphic algebras, hence 
the Hilbert series 
of $\twotwoqgenmat$ is 
\[\frac{(1-s)(1-t)(1-st)+st}
{(1-s)^2(1-t)^2(1-st)}.\]  
\end{theorem} 

\begin{proof} 
The elements $a,b,c,d,e$ are obviously central in $\twotwotracering$, 
the Cayley-Hamilton Theorem asserts  
$x^2=ax-c$, $y^2=by-d$, 
and its bilinearization says  
$yx=-xy+ay+bx-e$. 
With these equalities at hand, direct computation verifies that 
the elements $x,y,a,b,c,d,e+(1-q^{-2})xy$ satisfy 
the relations of $X,Y,A,B,C,D,E$. 
Therefore there exists an algebra homomorphism 
$\phi:\twotwoqtracering\to\twotwotracering$ 
mapping the generators to the prescribed images. These images clearly 
generate the whole $\twotwotracering$, hence 
$\phi$ is surjective. On the other hand, $\phi$ 
is a grading preserving map between two spaces with the same Hilbert series. 
So $\phi$ must be an isomorphism. 
Finally, note that under this isomorphism, $\twotwoqgenmat$ is mapped onto 
$\twotwogenmat$. The Hilbert series of $\twotwogenmat$ was determined in 
\cite{fhl}. 
\end{proof} 

\begin{remark} 
{\rm Note that the above isomorphism $\twotwoqtracering\to\twotwotracering$ 
is non-trivial in the sense that the subalgebra 
$\twotwoqmatinv I$ is not mapped into $\twotwomatinv I$: the element 
$E$ is mapped outside $\twotwomatinv I$. }
\end{remark} 

\begin{remark} 
{\rm The algebras $\genmat$ are highly non-trivial objects 
from the combinatorial 
point of view. For example, they are not finitely presented algebras: 
the ideal of relations among the generators is not finitely generated 
(although it is finitely generated as a T-ideal);  
see the papers \cite{dr}, \cite{f1}, \cite{pr2}, \cite{llb} 
for some computations 
in the $2\times 2$ case. 
Therefore the isomorphism 
$\twotwoqgenmat\cong\twotwogenmat$ looks a bit  surprising to us. }
\end{remark} 

Not all $\mcoact$-coinvariants are central in $\qsevmat$, 
and the algebra $\qmatinv$ is not commutative in general,
as the following Proposition shows.

\begin{proposition} \label{prop:non-comm} 
Suppose that $q^2\neq 1$. 
\begin{itemize}
\item[(i)] The elements $\qtr(XY)$ and $\qtr(XZ)$ 
do not commute in $\twothreeqsevmat$. 
\item[(ii)] The element $\qtr(XY)$ is not central in $\twotwoqsevmat$. 
\end{itemize} 
\end{proposition} 

\begin{proof} 
(i) The expansion of the commutator $[\qtr(XY),\qtr(XZ)]$ in the basis 
\eqref{eq:qsevmat22-basis}  is 
$\sum_{i,j,k,l}f^{ik}_{jl}y^i_jz^k_l$, 
where $f^{ik}_{jl}$ is a homogeneous quadratic expression of the variables 
$x^s_t$. 
We claim that $f^{22}_{11}$, the coefficient of $y^2_1z^2_1$, is 
non-zero. An inspection of \eqref{eq:2rea}, \eqref{eq:cross} shows 
that $y^2_1$ (respectively, $z^2_1$) appears on the right hand side of a 
relation if and only if $y^2_1$ (respectively, $z^2_1$) 
appears on the left hand side of the relation. 
Therefore the coefficient of $y^2_1z^2_1$ is the same as in 
$[qx^1_2y^2_1,qx^1_2z^2_1]$, hence $f^{22}_{11}=(1-q^2)x^1_2x^1_2$. 

(ii) This follows from the relations in Theorem~\ref{th:presentation} (iii), 
but it is easy to check directly, 
similarly to (i).  Write 
$[\qtr(XY),x^1_1]$ as 
$\sum f^i_j y^i_j$, where 
$f^i_j$ is a linear combination of the entries of $X$. 
We claim that $f^2_1\neq 0$. 
By the same observation on the relations we used in (i), 
we see that $f^2_1$ is the same as the coefficient of $y^2_1$ in 
$[qx^1_2y^2_1,x^1_1]$, so 
$f^2_1=q[x^1_2,x^1_1]\neq 0$. 
\end{proof} 

\begin{remark} 
{\rm If we endow the coordinate ring of pairs of $2\times 2$ matrices 
with the Poisson bracket that is the classical limit of the tensor product 
of two copies of the quantum coordinate ring of $2\times 2$ matrices, 
then a calculation similar to the above proof yields that 
$\twotwomatinv$ is not a Poisson subalgebra. 
So roughly speaking, there is no subalgebra in the tensor product of two
copies of the quantum coordinate ring of $2\times 2$ matrices that would lie 
above $\twotwomatinv$. This indicates that we are forced to move to
the reflection equation algebra and braided tensor products to find quantum
analogues of the rings of matrix invariants. }
\end{remark}

We continue with the ring theoretic study of our algebras. 

\begin{theorem}\label{th:2x2noether} 
\begin{itemize} 
\item[(i)] The algebra $\twomqsevmat$ is noetherian. 
\item[(ii)] Suppose that $q$ is not a root of unity. 
Then $\twomqmatinv$ is noetherian and is finitely generated as a 
$\mc$-algebra. 
\item[(iii)] Suppose that $q$ is not a root of unity. 
Then $\twomqtracering$ is finitely generated as a left (respectively right) 
$\twomqmatinv$-module. In particular, it is noetherian and is finitely
generated as a $\mc$-algebra. 
\end{itemize}
\end{theorem} 

\begin{proof} (i) One may build up $\twomqsevmat$ by adjoining step-by-step the 
variables $x(k)^i_j$ in the following order: 
\[x(1)^1_2,x(1)^2_2,x(1)^1_1,x(1)^2_1, 
x(2)^1_2,x(2)^2_2,x(2)^1_1,x(2)^2_1,\ldots.\] 
For notational simplicity, rename the members of this ordered sequence as 
$x_1,x_2,\ldots,x_{4m}$. 
Denote by $A_i$ the subalgebra generated by $x_1,\ldots,x_i$. 
Then it is easy to see from the defining relations of $\twomqsevmat$ that 
for all $i$, we have that $A_ix_{i+1}$ is contained in $A_i+x_{i+1}A_i$, 
and $x_{i+1}A_i$ is contained in $A_i+A_ix_{i+1}$. 
This implies successively that the algebras $A_i$ are noetherian 
(see 1.2.10 in \cite{mr}), and thus $\twomqsevmat$ is noetherian. 

(ii) The Haar functional of the cosemisimple Hopf algebra $\qsl$ can be 
used to construct a projection (Reynolds operator) 
$\qsevmat\to\qmatinv$ with certain good 
properties (see for example \cite{dl2} for details). 
Then both statements can be derived from (i) by the well known argument of 
Hilbert. 

(iii) 
In order to show that 
$\twomqtracering$ is a finitely generated left (right) module 
over $\twomqmatinv$,  
we recall from Remark~\ref{rem:1} the 
$\mc$-vector space isomorphism 
$\twomqtracering\cong (\twomqsevmat\otimes\xijspan^*)^{\twoqsl}$. 
This is clearly an isomorphism of left (right) $\twomqmatinv$-modules  
(note that $\twomqsevmat\otimes\xijspan^*$ is naturally a left and right 
$\twomqsevmat$-module). 
Now $\twomqsevmat\otimes\xijspan^*$ is a finitely generated module 
over $\twomqsevmat$, hence is noetherian by (i). Hence 
$(\twomqsevmat\otimes\xijspan^*)^{\twoqsl}$ 
is a finitely generated module over 
$\twomqsevmat^{\twoqsl}= \twomqmatinv$  
by a well known argument using the Reynolds operator 
$\twomqsevmat\to\twomqmatinv$ again; see for example page 71 in 
\cite{dol} for this argument. 
\end{proof}


\noindent M. Domokos: 
R\'enyi Institute of Mathematics, \\
Hungarian Academy of Sciences,\\ 
P.O. Box 127, 1364 Budapest, Hungary\\
E-mail: domokos@renyi.hu\\
\\
\noindent T. H. Lenagan: 
School of Mathematics, University of Edinburgh,
\\ James Clerk Maxwell Building, King's Buildings, Mayfield Road, 
\\Edinburgh EH9 3JZ, Scotland
\\E-mail: tom@maths.ed.ac.uk

\end{document}